\documentstyle[12pt,righttag,ctagsplt]{amsart}
\newtheorem{theorem}{Theorem}[section]

\newtheorem{prop}[theorem]{Proposition}

\theoremstyle{definition}
\newtheorem{defin}{Definition}
\newtheorem{example}[theorem]{Example}

\theoremstyle{remark}

\newtheorem{ack}{Acknowledgments} 

\numberwithin{equation}{section}

\newcommand{\bbR}{{\mathbb{R}}}
\newcommand{\bbC}{{\mathbb{C}}}

\newcommand{\bbK}{{\mathbb{K}}}

\newcommand{\calF}{{\cal{F}}}

\newcommand{\de}{\delta}
\newcommand{\e}{\varepsilon}

\newcommand{\f}{\varphi}

\renewcommand{\span}{\operatorname{span}}

\newcommand{\conv}{\operatorname{conv}}

\newcommand{\Ker}{\operatorname{Ker}}

\newcommand{\lb}{\label}

\newcommand{\lra}{\longrightarrow}

\newcommand{\wtw}{if and only if }

\newcommand{\ONTO}{\buildrel {\mbox{\small onto}}\over \longrightarrow}

\makeatletter
\def\@currentlabel{2.1}\label{e:dispaa}
\def\@currentlabel{2.21}\label{e:dispau}
\def\@currentlabel{2.22}\label{e:dispav}
 \def\@currentlabel{2.23}\label{e:dispaw}
\def\@currentlabel{2.24}\label{e:dispax}

\makeatother

\makeatletter
\def\alphenumi{%
  \def\theenumi{\alph{enumi}}%
  \def\p@enumi{\theenumi}%
  \def\labelenumi{(\@alph\c@enumi)}}
\makeatother
\newsymbol\subsetneq 2328
\newsymbol\subsetneqq 2324
\newsymbol\supsetneq 2329
\newsymbol\nsubseteq 232A

\begin{document}
\title{A note on Banach--Mazur problem}

\author{Beata Randrianantoanina$^*$}\thanks{$^*$Participant, NSF Workshop
in Linear Analysis and Probability, Texas A\&M University}

\address{Department of Mathematics and Statistics
\\ Miami University \\Oxford, OH 45056}

 \email{randrib@@muohio.edu}


\subjclass{46C15,46B04,46B20}

\begin{abstract}
We prove that if $X$ is a real Banach space, with $\dim X\geq 3$, which
contains a subspace of codimension 1 which is 1-complemented in $X$ and
whose group of isometries is almost transitive then $X$ is isometric to a
Hilbert space.
This partially answers the Banach-Mazur rotation problem and generalizes
some recent related results.
\end{abstract}
\maketitle

\section{Introduction}

In 1930's  Banach and Mazur
\cite{Ban} (see also \cite[Problem~9.6.2]{Rol})
 posed a   problem whether every separable Banach space with 
transitive group of
isometries has to be isometric to a Hilbert space.  
Here we say that the
group of isometries of a Banach space $X$ is {\it transitive} if for every
$x,y\in X$ with $\|x\|=\|y\|=1$ there exists an isometry 
$T:X\ONTO X$ such
that
$$
Tx=y.
$$

Mazur \cite{M38} answered this problem positively for finite
dimensional Banach spaces $X$ and Pe\l czy\'nski and Rolewicz \cite{PR62}
showed
that the answer is negative  when $X$ is not assumed to be separable.
The case of infinite dimensional separable spaces remains open despite
active research in the area, see \cite{Rol} and a recent survey \cite{Cab97}.

Closely related to the notion of transitivity are the notions of 
almost transitivity
and convex transitivity.
We say that a group of isometries of $X$ is   {\it almost transitive}
(resp. {\it convex transitive}) if for every $x$ in the unit sphere of $X$, 
$S_X=\{x:\|x\|=1\}$ the orbit of $x$ i.e. the
set $G_x=\{Tx:T$ isometry of $X\}$ is dense in $S_X$ 
(resp. $\conv(G_x)$ is
dense in $S_X$).
Sometimes we will abuse language and say that a space 
$X$ is   {almost transitive}
(resp. {convex transitive} or transitive) provided 
the group of isometries of $X$ is    { almost transitive}
(resp. {convex transitive} or transitive).

Spaces with almost transitive and convex transitive groups of isometries
have been actively studied. It is known that 
there exist 
non-Hilbertian separable Banach spaces with 
almost transitive   groups of isometries, for example 
$L_p[0,1], \ 1\le p<\infty,$
 are such spaces \cite{PR62}, see also \cite{Rol} and see \cite{GJK} for 
detailed study which function spaces are almost transitive.
Several questions have been posed to find additional conditions on a 
Banach space
$X$ which together with almost transitivity, or with just convex
  transitivity, imply that $X$ is isometric to a Hilbert space. 
Maybe the most famous conjecture of this type is the conjecture of Wood
\cite{Wood82} that if $C_0(L)$ is almost transitive in its natural  supremum
norm then $L$ is a singleton, i.e. $C_0(L)$ is one dimensional. 
Wood's conjecture
is still open despite recent active research in the area, see \cite{Cab97}.

The main theorem of the present paper is the following:

\begin{theorem} \lb{intro} (see Theorem~\ref{main} below)
Suppose that $X$ is a real Banach space, with dim$X\geq 3$,
which contains a 1-complemented
hyperplane  and whose group of isometries is
almost transitive.  Then $X$ is isometric to a  Hilbert space.
\end{theorem}

Our method of proof relies on the theory of numerical ranges \cite{BD1,BD2}.
We postpone the proof to the next section and now we will discuss the
connections with existing results in the literature.

Note that spaces $L_p[0,1],\ \ 1\le p<\infty$, do not have 1-complemented
hyperplanes (see e.g. \cite{KR,pams}) and, as mentioned above, they are 
almost transitive.

Theorem~\ref{intro} generalizes a recent result of Skorik and Zaidenberg:
 
\begin{theorem} \cite{SZ97} \lb{SZ}  
Suppose that $X$ is a real Banach space which contains an isometric
reflection and whose group of isometries is almost transitive.  Then $X$ is
isometric to a Hilbert space.
\end{theorem}

Here we say  that an operator $T$ is {\it a reflection on $X$} 
if there exist $e\in X, \ e^*\in X^*$ with
$e^*(e)=1$ so that 
$$T(x)= x-2e^*(x)e.$$ 
If this happens we write $T=S_{e,e^*}$, and if $ S_{e,e^*}$ is an isometry
we say that $e$ is an {\it isometric
reflection vector} in $X$.

To see that Theorem ~\ref{intro} is more general than Theorem~\ref{SZ}
we observe that if a space $X$ admits
an isometric reflection operator then $X$ contains a 1-complemented
hyperplane but not vice-versa. 

Indeed,  assume that   for some $e\in S_X$ there exists $e^*\in X^*$ 
with $e^*(e)=1$ so
that an operator $S_{e,e^*}: x\mapsto x-2e^*(x)e$ is an isometric
reflection in $X$. Then for all $x\in X$ we have:

$$\|x-e^*(x)e\|=\|\frac 12 ((x-2e^*(x)e)+x)\|\leq \frac 12
(\|x-2e^*(x)e\|+\|x\|)=\|x\|
$$
Thus $x\mapsto x-e^*(x)e$ is a contractive projection in $X$ onto the hyperplane
$\Ker(e^*)\subset X$. 

On the other hand consider the
two dimensional real space $X$ whose unit sphere is the convex hull
of the points (1,0), (1/2,1), (-1,1), (-1,0), (-1/2,-1), (-1,-1) as 
sketched below.
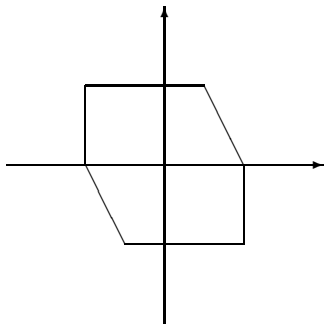
\vspace{50mm}

\hspace{50mm}
\begin{picture}(90,90)(-75,-75)
\put(-60,0){\vector(1,0){120}}
\put(0,-60){\vector(0,1){120}}
\put(30,0){\line(-1,2){15}}
\put(15,30){\line(-1,0){45}}
\put(-30,30){\line(0,-1){30}}
\put(-30,0){\line(1,-2){15}}
\put(-15,-30){\line(1,0){45}}
\put(30,-30){\line(0,1){30}}
\end{picture}

Then $X$ contains no isometric reflections but clearly every hyperplane is 
1-complemented in $X$.

It is not difficult to construct spaces of arbitrary dimension which
contain no isometric reflections but which do contain 
 1-complemented subspaces of
codimension 1.

Statements similar to Theorem~\ref{SZ} have been recently studied by 
J. Becerra Guerrero, F. Cabello Sanchez and A.
Rodriguez Palacios.
F. Cabello Sanchez \cite{Cab96}  showed that Theorem~\ref{SZ} is
 valid in
complex Banach spaces.  J. Becerra Guerrero  and A.
Rodriguez Palacios linked this result with the following characterization of
Hilbert spaces due to Berkson \cite{Berk72} and Kalton and Wood \cite{KW}:

\begin{theorem}\lb{BKW}
Let $X$ be a complex Banach space. If $x\in X$ is such that $\span\{x\}$ is
the range of a hermitian projection in $X$ then $x$ is called a hermitian element in
$X$.

If every nonzero element of $X$ is hermitian in $X$ then $X$ is a Hilbert
space.
\end{theorem}

J. Becerra Guerrero  and A.
Rodriguez Palacios \cite{BR} observed that  an element  $x\in X$ is hermitian in $X$
\wtw $x$ is an isometric reflection vector   in $X$. Thus in the 
complex case we have the following stronger version of Theorem~\ref{SZ}:

\begin{theorem} \cite[Theorem~6.4]{KW} \lb{convex}
Let $X$ be a complex Banach space. If $ X$ is convex transitive and $X$
contains an isometric reflection vector then $X$ is a Hilbert
space.
\end{theorem}

J. Becerra Guerrero  and A.
Rodriguez Palacios generalized  
Theorem~\ref{SZ} as follows:

\begin{theorem} \cite{BR99b} \lb{BR}
Let $X$ be a real or complex Banach space. If there exists a nonrare set in $S_X$
consisting of
 isometric reflection vectors then $X$ is a Hilbert
space.
\end{theorem}

We do not know whether Theorem~\ref{intro} can be generealized for 
convex transitive spaces (to obtain an analogue of Theorem~\ref{convex}).
However the following example illustrates that an analogue of
Theorem~\ref{BR} for norm-one complemented hyperplanes fails in a very
strong way. Namely we have:

\begin{example} \lb{epsilon}
For every $\e>0$ there exists a 3-dimensional Banach space $X_\e$ which is
not isometric to a Hilbert space and such that the set $F$ of functionals
$f\in S_{X_\e^*}$ with $\Ker f$ norm-one complemented in $X_\e$ is open
in $S_{X_\e^*}$ and $\mu(S_{X_\e^*}\setminus F)<\e,$ (here $\mu$ denotes
the Lebesgue measure on $S_{X_\e^*}\subset \bbK^3, \bbK= \bbR $ or $\bbC$).

Moreover $X_\e$ can be chosen to be uniformly convex.
\end{example}

\begin{pf}
Let $\e>0$. Fix $\de>0$ so that 
$$\mu(\{(x_1,x_2,x_3)\in S_{\ell_2}\ :\ |x_3|\ge 1-\de\})<\e.$$

Let $\f:(\bbR_+)^3\lra \bbR$ be a convex continuous function such
that
\begin{equation*}
\f(t_1,t_2,t_3)\begin{cases}
=1 & {\text{if}}\  (t_1,t_2,t_3)=(0,0,1)\\
=\sqrt{t_1^2+t_2^2+t_3^2} & {\text{if}}\  t_3\le 
(1-\de)\sqrt{t_1^2+t_2^2+t_3^2} \\
>\sqrt{t_1^2+t_2^2+t_3^2} & {\text{if}}\  t_3> 
(1-\de)\sqrt{t_1^2+t_2^2+t_3^2}. 
\end{cases}
\end{equation*}
We can additionally require $\f$ to have any desired degree of
smoothness.

We define a norm on $\bbK^3$ using function $\f$:
$$\|(x_1,x_2,x_3)\|_{X_\e}=\f(|x_1|,|x_2|,|x_3|).$$
Then $S_{X_\e}\subsetneq S_{\ell_2}$ and
$S_{X_\e^*}\cap S_{\ell_2}\supset 
\{(x_1,x_2,x_3)\in S_{\ell_2}\ :\ |x_3|\le 1-\de\}$
as illustrated in the figure below.

\vspace{9mm}

\begin{quote}
{\it Here should come another figure which doesn't run well
in TeX so it is not included here.
 The ps file of this paper
which includes all figures is available at}\hspace{3mm}
{\tt http://www.users.muohio.edu/randrib/bm3.ps}
\end{quote}

\vspace{9mm}

Thus, if $f\in S_{X_\e^*}\cap S_{\ell_2}$ then
$\Ker f \cap S_{X_\e}\subset S_{\ell_2}$ and the orthogonal projection
$P$ onto $\Ker f$ has norm 1 in $X_\e$. Hence
$$F=\{f\in S_{X_\e^*}\ : \ \Ker f {\text{ is 1-complemented in }} X_\e\}
\supset S_{X_\e^*}\cap S_{\ell_2}.$$

Thus $\mu(S_{X_\e^*}\setminus F)<\e,$ as desired. 
\end{pf}

Theorem~\ref{intro} relies on the real
analogue of the notion of hermitian elements
(see Definitions~1 and 2 below) whose existence in $X$ is equivalent
to the existence of norm one complemented hyperplanes in $X$.
We note here that our proof is much shorter than the existing proofs of
Theorem~\ref{SZ}.

\section{Proofs of main results}

We begin with definitions of real analogues of hermitian elements which
 were introduced by Kalton and the
author \cite{KR} based on ideas of P. H. Flinn, cf. \cite{R84}.

\begin{defin} \cite{R84} We say that an operator $T:X\to X$ is {\it
numerically positive} if for all $x\in X$ there exists a $x^*\in X^*$, or,
equivalently, for all $x^*\in X^*$ with $\|x^*\|^2=\|x\|^2=x^*(x)$ and
$x^*(Tx)\geq 0$, i.e. the numerical range of $T$ is contained in 
${\Bbb{R}}_+\cup\{ 0 \}$ (cf. \cite{BD1,L61}). \end{defin}

\begin{defin} \cite{KR} (based on ideas of P.H. Flinn \cite{R84})  We say
that $u\in X$ is {\it a Flinn element} if there exists a numerically
positive projection $P:X\ONTO \span\{u\}$ i.e. if there
exists $f\in X^*$ with $f(u)=1$ and such that the map $f\otimes u$, defined
by $x\mapsto f(x)u$, is numerically positive.  We say then that $(u,f)$ is
a {\it Flinn pair}.\end{defin}

The set of all Flinn elements of $X$ will be denoted by $\calF (X)$.

We list few straightforward
properties of Flinn elements which are important for the future use.  We include their
short proofs for completeness.

\begin{prop}   \cite[Lemma~1.4]{R84} \lb{contractive} A projection
$P:X\lra X$, $P\ne I$, is numerically positive if and only if $\|I-P\|=1$.
\end{prop}

\begin{pf} 
If $\|I-P\|=1$ and $x\in X$, $x^*\in X^*$ are such that
 $1=\|x^*\|^2=\|x\|^2=x^*(x)$, then
$$1\ge x^*((I-P)(x))= 1-x^*(P(x)).$$
Thus $x^*(P(x))\ge 0$ and $P$ is numerically positive.

To see the implication in the other direction we rely on the result of 
Lumer and Phillips \cite{LP61} that operator $P$ is numerically positive
\wtw $\|\exp(-tP)\|\le 1$ for all real $t\ge 0$. We have for 
all real $t\ge 0$:
$$\exp(-tP)=I+\sum_{j=1}^\infty\frac{(-t)^jP^j}{j!} = 
I+\left(\sum_{j=1}^\infty\frac{(-t)^j}{j!}\right)P= I+(e^{-t}-1)P,$$
where the second equality holds because $P$ is a projection.
Thus by result of Lumer and Phillips 
if $P$ is numerically positive we have
$$\|I-P\|=\lim_{t\to\infty}\|\exp(-tP)\|\le1.$$
Since $P\ne I$ and $I-P$ is  a projection we get $\|I-P\|=1$.
\end{pf}

\begin{prop} \lb{isoflinn} \cite[Proposition~3.2]{KR}
Suppose that $T:X\lra X$ is a surjective isometry.  Then $T(\calF (X))=\calF
(X)$. \end{prop}

\begin{pf} If $u\in\calF (X)\setminus \{0\}$ then there exists $f\in X^*$
such that the projection
$P:X\lra\span\{u\}$ defined by $x\mapsto f(x)u$ is numerically positive.
Then the projection $Q:x\mapsto ((T^*)^{-1} f)
(x)\cdot T(u)$ establishes the fact that $T(u)$ is Flinn.
\end{pf}

\begin{prop} \lb{closed} \cite[Proposition~3.1]{KR} The set $\calF (X)$ is
closed. \end{prop}

\begin{pf} Suppose $u_n\in \calF (X)$ with lim$\|u_n-u\|=0$.  It suffices
to consider the case when $\|u_n\|\neq 0$ and $\|u\|\neq 0$.  Then there
exist $f_n\in X^*$ so that $f_n\otimes u_n$ is a numerically positive
projection.  Thus $\|f_n\otimes u_n\|=\|f_n\|\ \|u_n\|\leq 2$.  Thus
$\|f_n\|\leq 2$sup$(1/ \|u_n\|)$.  By Alaoglu's theorem $(f_n)$
has a weak$^*$-cluster point $f$ and clearly $(u,f)$ is a Flinn pair.
\end{pf}

\begin{prop}\lb{subspace} Suppose $u\in \calF (X)$ and that $Y$ is a subspace of $X$
such that $u\in Y$.  Then $u\in \calF (Y)$.\end{prop}

\begin{pf} Without loss of generality $u\neq 0$.  Since $u\in \calF (X)$
there exists $f\in X^*$ with $f(u)=1$ and such that the map $x\mapsto
f(x)u$ is numerically positive.  Consider $g=f|_Y\in Y^*$ and the map
$Q:Y\lra$span$\{u\}\subset Y$ defined by
$$
Q(y)=g(y)u.
$$
By Hahn-Banach Theorem for every $y\in Y$ and every $y^*\in Y^*$ with
$\|y^*\|^2=\|y\|^2=y^*(y)$ there exists $\tilde y^*\in X^*$ with $\|\tilde
y^*\|=\|y^*\|=\|y\|$ and $\tilde y^*|_Y=y^*$.  Thus we get, since $(u,f)$
is a Flinn pair in $X$:
$$
g(y)y^*(u)=f(y)\tilde y^*(u)\geq 0
$$
and $g(u)=f(u)=1$.  Hence $(u,g)$ is a Flinn pair in $Y$.
\end{pf}

Now we are ready for the real analogue of Theorem~\ref{BKW} 
(\cite[Theorem 2.22]{Berk72},
\cite[Corollary 4.4]{KW}).

\begin{theorem} \lb{flinn}
Suppose that $X$ is a real Banach space $\dim  X\geq 3$ and $\calF (X)=X.$
Then $X$ is isometric to a Hilbert space.
\end{theorem}

\begin{pf} Since Hilbert spaces are characterized by the parallelogram
identity   it is enough to show that the result holds for all
3-dimensional subspaces of $X$ (cf. \cite[(1.4')]{Amir}).  Suppose that
$Y\subset X$ is a real Banach space with dim$Y=3$.  Then, by 
Proposition~\ref{subspace},
$\calF (Y)=Y$ i.e. for every $u\in Y$ there exists an $f\in Y^*$ such that
$f(u)=1$ and the map $f\otimes u$ is a numerically positive projection in
Y, and also in $Y^*$.

Hence, by Proposition~\ref{contractive}, the map $I-f\otimes u$ is a contractive projection
of $Y^*$ onto $\Ker u\subset Y^*$.  Since $\dim Y=3$, we conclude that
every  2-dimensional subspace of $Y^*$ is contractively complemented.
Thus $Y^*$ is isometric to a Hilbert space by the following criterion due to
Kakutani:

\begin{theorem}  \cite{Ka39} (see also \cite[p.~99]{Amir}) Suppose that $Z$
is a Banach space of dimension at least 3.  Then $Z$ is isometric to a
Hilbert space if and only if for every 2-dimensional subspace $F$ of $Z$
there is a norm one linear projection $P:Z\lra F$.\end{theorem}

Therefore $Y$ is isometric to Hilbert space and the proof is finished.
\end{pf}

Next is our main theorem.

\begin{theorem}\lb{main} 
Suppose that $X$ is a real Banach space, with dim$X\geq 3$,
which contains a 1-complemented
hyperplane  and whose group of isometries is
almost transitive.  Then $X$ is isometric to a  Hilbert space.
\end{theorem}

\begin{pf}
When $X$ contains a subspace of codimension 1 which is
1-complemented in $X$ then, by Proposition~\ref{contractive}, $X$ contains a nonzero Flinn
element. Let
 $u\in S_X$ be a Flinn element.  By 
Proposition~\ref{isoflinn}
$$
G_u=\{Tu:T \text{ isometry of } X \text{ onto } X\}\subset \calF (X)
$$
and by almost transitivity of $X,\ G_u$ is dense in $S_X$.  Since, by
Proposition~\ref{closed}, $\calF(X)$ is closed we obtain:
$$
\calF (X)\cap S_X=S_X.
$$

Thus, $\calF (X)=X$ and the result follows by Theorem~\ref{flinn}.
\end{pf}

\begin{ack} I wish to thank Professor G. Wood for telling me about the
current status of work on problems related to Banach-Mazur problem, and to
Professors F. Cabello Sanchez and A. Rodriguez Palacios for providing me
with copies of their preprints.
\end{ack}

\noindent
{\it Added in proof.} After this paper has been completed and accepted for
publication P. L. Papini has pointed out to me the reference \cite{O82} which
contains a result analogous to Theorem~1.1 with the additional assumptions that
$X$ is reflexive and transitive.


\end{document}